\documentstyle[12pt,fullpage]{article}

\input{psfig}
\def\IMSmarkvadjust{0 pt}
\def\IMSmarkhadjust{0 pt}
\def\IMSmarkhpadding{0 pt}
\def\IMSpubltext{Published in modified form:}
\def\SBIMSMark#1#2#3{
 \font\SBF=cmss10 at 10 true pt
 \font\SBI=cmssi10 at 10 true pt
 \setbox0=\hbox{\SBF \hbox to \IMSmarkhpadding{\relax}
                Stony Brook IMS Preprint \##1}
 \setbox2=\hbox to \wd0{\hfil \SBI #2}
 \setbox4=\hbox to \wd0{\hfil \SBI #3}
 \setbox6=\hbox to \wd0{\hss
             \vbox{\hsize=\wd0 \parskip=0pt \baselineskip=10 true pt
                   \copy0 \break%
                   \copy2 \break% 
                   \copy4 \break}}
 \dimen0=\ht6   \advance\dimen0 by \vsize \advance\dimen0 by 8 true pt
                \advance\dimen0 by -\pagetotal
	        \advance\dimen0 by \IMSmarkvadjust
 \dimen2=\hsize \advance\dimen2 by .25 true in
	        \advance\dimen2 by \IMSmarkhadjust

  \openin2=publishd.tex
  \ifeof2\setbox0=\hbox to 0pt{}
  \else 
     \setbox0=\hbox to 3.1 true in{
                \vbox to \ht6{\hsize=3 true in \parskip=0pt  \noindent  
                {\SBI \IMSpubltext}\hfil\break
                {\bf Markov partitions and Feigenbaum-like mappings}. 
Commun. Math. Phys. {\bf 171} (1995), pp. 351-363.  
                \vfill}}
  \fi
  \closein2
  \ht0=0pt \dp0=0pt
 \ht6=0pt \dp6=0pt
 \setbox8=\vbox to \dimen0{\vfill \hbox to \dimen2{\copy0 \hss \copy6}}
 \ht8=0pt \dp8=0pt \wd8=0pt
 \copy8
 \message{*** Stony Brook IMS Preprint #1, #2. #3 ***}
}

\begin{document}
\def\IMSmarkvadjust{-10pt}
\SBIMSMark{1991/19a}{June 1991}{}
\thispagestyle{empty}
%\vspace*{1in}
\centerline {\bf On the quasisymmetrical classification of infinitely
renormalizable maps}

\vskip10pt
\centerline {\sc I. Maps with Feigenbaum's topology}

\vskip15pt
\centerline{by Yunping Jiang }

\vskip10pt
\centerline{June, 1991}

\vskip30pt
\noindent {\bf \S 0 Introduction}

\vskip5pt
We begin by considering the set of infinitely renormalizable unimodal maps 
on the interval $[-1, 1]$. A function $f$
defined on $[-1,1]$ is said to be {\em unimodal} if it is continuous,
increasing on $[-1, 0]$, decreasing on $[0,1]$ and  
symmetric about $0$,  and if it fixes $-1$ and maps $1$ to $-1$ .  Moreover,
it is said to be {\em renormalizable} if there is an integer $n>1$
and a subinterval $I\neq [-1, 1]$ containing $0$ such that $f^{\circ n}(I)
\subseteq I$. We will assume that one endpoint $q$ of $I$ is fixed by $f^{\circ n}$
and $I$ is symmetric about $0$. If we normalize $I$ to $[-1, 1]$ by the linear
map $\alpha$, which maps $q$ to $-1$, then ${\cal R}(f) = \alpha \circ
f^{\circ n} \circ \alpha^{-1}$ is unimodal, too. This map ${\cal R}$ will be 
called the {\em renormalization operator}. To fix our notation, we will assume that $n$ is the minimum
such integer and call it the {\em return time}.

\vskip5pt
A unimodal map $f$ is {\em infinitely renormalizable} if every 
${\cal R}^{\circ k}(f)$ is renormalizable, say with {return time} $n_{k}$. 
Furthermore, $f$ is of {\em bounded type} if 
all the return times are less than a constant
integer, otherwise, $f$ is of {\em unbounded type}. 
In particular, we call $f$ a Feigenbaum map if
all the return times are $2$. 
A well-known example $f=q_{\lambda_{\infty}}$ of a Feigenbaum map 
\cite{ct,f,vsk} is obtained by
period-doubling cascade in the family $\{ q_{\lambda }(z)
= -(1+\lambda)z^{2} +\lambda \}_{0\leq \lambda \leq 1}$. 

\vskip5pt
Let ${\cal U}$ be the space of unimodal maps $f=h\circ Q_{t}$ where
$Q_{t}(x)=-|x|^{t}$ for some $t>1$ and $h$ is a diffeomorphism (a unimodal map
can always be written in this form by some smooth change of coordinate \cite{j3}). We may assume
either $h$ is a $C^{3}$-diffeomorphism with nonpositive Schwarzian 
derivative \cite{cec} or
$h$ is a $C^{1+Z}$-diffeomorphism \cite{su}. However, the
smoothness is not important in this paper as long as $h$ satisfies the
distortion properties discussed in \cite{su}. To avoid many technical
notations, {\em henceforth, we will assume that $h$ is a
$C^{3}$-diffeomorphism
with
nonpositive Schwarzian derivative.} We note that the Schwarzian
derivative $S(h)$ of a $C^{3}$-diffeomorphism is 
\[ S(h) = \frac{ h'''}{h'}-\frac{3}{2}\Big( \frac{
h''}{h'}\Big)^{2}.\]

\vskip5pt
Suppose $f$
and $g$ in ${\cal U}$ are both infinitely renormalizable and topologically conjugate. 
Suppose $H$ from $[-1,1]$ to $[-1, 1]$ is the conjugacy between $f$ and $g$, that is, 
$f\circ H =H\circ g$. The map  
$H$ is said to be quasisymmetric \cite{al} if there is a constant $M>0$ such that 
\[ M^{-1}\leq \frac{|H(x)-H(z)|}{|H(z)-H(y)|} \leq M\]
whenever $x$ and $y$ are in $[-1,1]$ and $z=(x+y)/2$ is the midpoint between
$x$ and $y$. 
 
\vskip10pt
{\sc Conjecture 1 (Sullivan).} 
{\em The homeomorphism $H$ is quasisymmetric.}

\vskip5pt
We would like to investigate this conjecture in the three cases, the
Feigenbaum case, the bounded case and the unbounded case.  In this paper we 
prove this conjecture for the Feigenbaum case as follows. We note that any
two Feigenbaum maps are topologically conjugate. 

\vskip5pt 
{\sc Theorem 1.} {\em Suppose $f$ and $g$ in ${\cal U}$ are two Feigenbaum
maps and $H$ is the conjugacy between $f$ and $g$. Then $H$ is
quasisymmetric.}

\vskip5pt
The idea to prove this theorem is that we construct a sequence of
nested partitions on $[-1, 1]$ by using all the periodic points of 
a Feigenbaum map and the  
preimages of all the periodic points under 
iterates of this map.  We show that this sequence of nested partitions 
has bounded geometry
and bounded nearby geometry properties \cite{j2}. The reader may compare the
construction of the sequence of nested partitions here with the construction
of the sequence of nested
partitions in \cite{y} (see also \cite{h}) for a nonrenormalizable quadratic
polynomial.      

\vskip5pt
The proof of Conjecture 1 for infinitely renormalizable maps of bounded type
should not have an essential difference from the proof in this paper. We will
do it in a short paper \cite{j4}.
However, a proof of Conjecture 1 for infinitely renormalizable maps 
of unbounded
type may have an essential difference from the proof in this paper. It is
still an interesting research problem. 

\vskip5pt
For complex quadratic-like maps \cite{dh} in ${\cal U}$ which are infinitely
renormalizable of bounded type, Conjecture 1 has been proven by Sullivan 
\cite{su} using a completely different method. I was also told that recently, 
Paluba \cite{p} reached another proof of Conjecture 1 for infinitely renormalizable
maps of bounded type. Jakobson and Swiatek \cite{s,js} showed some
other interesting results about quasisymmetric classification of unimodal 
maps.

\vskip20pt
{\bf Acknowledgement.} The author would like to thank D. Sullivan for
introducing him this problem and to thank M. Shishikura for
introducing him the recent work of Yoccoz \cite{y} (see also \cite{h}) which inspired
the work in this paper. He would also like to thank J. Milnor for reading
and correcting this manuscript.

\vskip30pt
\noindent {\bf \S 1 The proof of Theorem 1}

\vskip5pt
We prove Theorem 1 by several lemmas.

\vskip5pt
Suppose $f=h\circ Q_{t}$, for some $t>1$, in ${\cal U}$ is a Feigenbaum map. 
Let $p_{0}=-1$ and $p_{1}\in (0, 1)$ be the fixed points of
$f$. Inductively, let $p_{n}\in (-p_{n-1}, 0)$ or $(0, -p_{n-1})$ be the
fixed point of $f^{\circ 2^{n-1}}$ for an integer $n>1$. Let $I_{n}$ be the interval bounded
by $p_{n}$ and $-p_{n}$ (see Figure 1) and $c_{n}=f^{\circ n}(0)$ 
is the $n^{th}$ critical value of $f$ for
$n=0$, $1$, $\cdots $.

\vskip5pt
\centerline{\psfig{figure=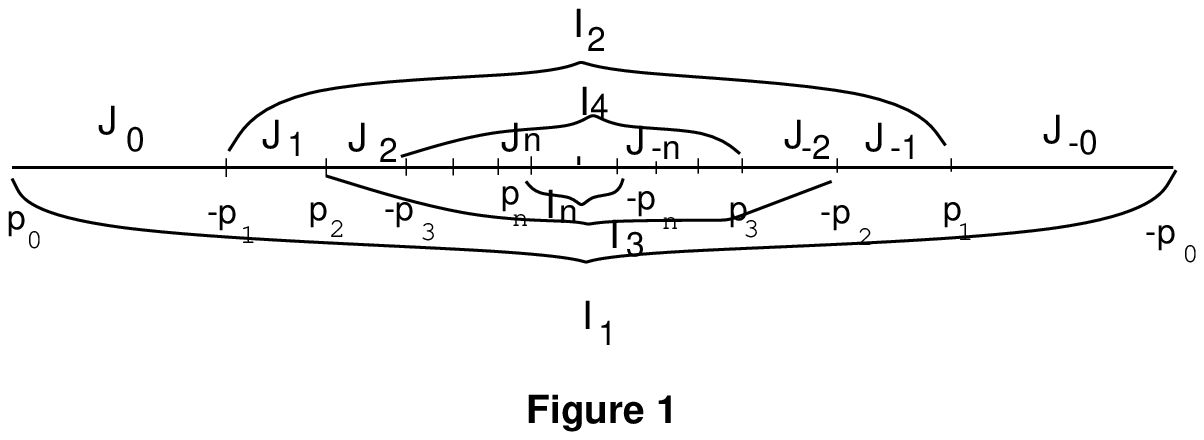}}

\vskip5pt
Suppose $L_{n}$ is the image of $I_{n-1}$ under $f^{\circ 2^{n-1}}$ and
$T_{n}$ is the interval bounded by the points $p_{n-1}$ and $p_{n}$. Let 
$M_{n}$ be the complement of $T_{n}$ in $L_{n}$ (see Figure 2). Then $M_{n}$
is the interval bounded by $p_{n}$ and $c_{2^{n-1}}$.

\vskip5pt
\centerline{\psfig{figure=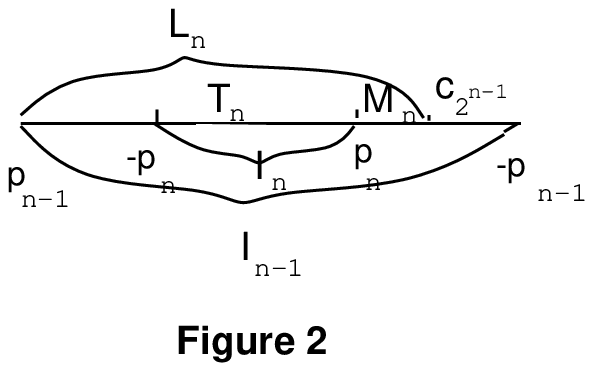}}

\vskip5pt 
{\sc Lemma 1.} {\em There is a constant $C_{1}=C_{1}(f)>0$ such that
$|M_{n}|/|I_{n}|\geq C_{1}$ for all the integers $n \geq 0$ .} 

\vskip5pt
{\sl Proof.} This lemma is actually proved in \cite{su}
by using the techniques
such as the smallest
interval and shuffle permutation on the intervals. 

\vskip5pt
{\sc Lemma 2.} {\em There is a constant $C_{2}=C_{2}(f)>0$ such that
$|I_{n}|/|I_{n-1}|\geq C_{2}$ for all the integers $n\geq 0$.}
 
\vskip5pt
We first prove a more general result as follows.
Suppose $h$ is a $C^{2}$-diffeomorphism. We use 
\[ N(h) = \frac{h''}{\Big( h' \Big)^{2}}\]
to denote the nonlinearity of $h$.
Let ${\cal K} = {\cal K}(t, K)$, for fixed numbers $t>1$ and $K>0$, 
be the subspace of maps $f=h\circ Q_{t}$ in ${\cal U}$ such that
\[ |\Big( N(h)\Big) (x)| \leq K\]
for all $x$ in $[-1,0]$. Remember that $Q_{t}(x)=-|x|^{t}$.

\vskip5pt
{\sc Lemma 3.} {\em There is a constant $C_{3}=C_{3}(t, K)>0$ such that
$f(0)=c_{1}(f) \geq C_{3}$ for all Feigenbaum maps $f$ in ${\cal K}$.}

\vskip5pt
{\sl Proof.} 
Suppose $f=h\circ Q_{t}$ is a map in ${\cal K}$. Since $h$ is a 
$C^{3}$-diffeomorphism with nonpositive
Schwarzian derivative, we can compare $h$ with some M\"obius transformations
$m(x)=(ax+b)/(cx+d)$. We note that all the M\"obius transformation 
have zero Schwarzian
derivatives. 

\vskip5pt
Let $m$ be the M\"obius transformation satisfying

\vskip3pt
(a) $m(-1) =f(-1)=-1$ and $m(0)=f(0)$, and 

\vskip3pt
(b) $\Big( N(m) \Big) (-1) =K$. 

\noindent Then 
\[ m(x) =\frac{(f(0) +1)(1+\frac{K}{2})x +f(0)}{\frac{K}{2}(f(0)+1)x+1}.\]

\vskip5pt
Since $m(-1)=h(-1)$, $m(0)=h(0)$ and $\Big( N(m)\Big) (-1) \geq \Big(
N(h)\Big) (-1)$, we have that $h(x)\geq m(x)$ for all $x$ in $[-1, 0]$. 

\vskip5pt
Suppose, in the moment, $w=f(0)>0$ in $m$ is a variable. Consider the
equation
\[  m\Big( -\Big(
m(0) \Big)^{t}
\Big)=\frac{(w+1)(1+\frac{K}{2})(-w^{t})+w}{\frac{K}{2}(w+1)(-w^{t})+1} =0. \]
Let $C_{3}=C_{3}(t,K)>0$ be the solution of this equation for $w$.

\vskip5pt
For any Feigenbaum map $f=h\circ Q_{t}$ in ${\cal K}$, 
we know, 
from the property of renormalization and $m(x)\leq h(x)$ for
all $x$ in $[-1,0]$, that $m\Big( -\Big( m(0)\Big)^{t}\Big) < 0$. It implies
that $f(0)=c_{1}(f)> C_{3}$. This proves Lemma 3.

\vskip5pt
Suppose $f=h\circ Q_{t}$, for some $t>1$, in ${\cal U}$ is a Feigenbaum map and $f_{n}= h_{n}\circ Q_{t}$ is the $n^{th}$ renormalization ${\cal R}^{\circ
n}(f)$ of $f$ for an integer $n\geq 0$. We note that the graph of $f_{n}$
is the rescale of the graph of the restriction of $f^{\circ 2^{n}}$ to $I_{n}$. 
\vskip5pt
{\sc Lemma 4}. {\em There is a constant $C_{4}=C_{4}(f)>0$ such that  
\[ |\Big( N(h_{n})\Big) (x)| \leq C_{4} \]
for all $x$ in $[-1, 0]$ and $n \geq 0$.}

\vskip5pt
{\sl Proof.} It is the a prior bound proved in \cite{su}

\vskip5pt
{\sl Proof of Lemma 2}. It is now a direct corollary of Lemma 1, 
Lemma 3 and Lemma 4 for $K=C_{4}$.

\vskip5pt
The set of the nested intervals $\{ I_{0}$, $I_{1}$, $\cdots $, $I_{n}$,
$\cdots \}$ gives a partition of $[-1,1]$ as follows.
Let $J_{0}$ and $J_{-0}$ be two connected components of $I_{0}\setminus
I_{1}$. Inductively, let $J_{n}$ and $J_{-n}$ be two connected components of
$I_{n}\setminus I_{n+1}$. Then the set $\eta_{0}= \{ J_{0}$, $J_{-0}$, $J_{1}$, $J_{-1}$,
$\cdots \}$ forms a partition of $I_{0}$ (see Figure 1). 

\vskip5pt
Now we are going to define a Markov map $F$ induced by $f$. Let $F$ be a
function of $[-1,1]$ defined as 
\[ F(x) =\left\{ \begin{array}{ll}
                         f(x) ,& x\in J_{0} \cup J_{-0},\\
                         f^{\circ 2}(x) ,& x\in J_{1}\cup J_{-1},\\
                         \vdots &  \\
                         f^{\circ 2^{n}}(x), & x\in J_{n}\cup J_{-n},\\
                         \vdots & .
                         \end{array}
                \right.\]

\vskip5pt
It is clear that $F$ is a Markov map in the sense that the image of every
$J_{k}$ under $F$ is the union of some intervals in the partition $\eta_{0}$
(see Figure 3).

\vskip5pt
\centerline{\psfig{figure=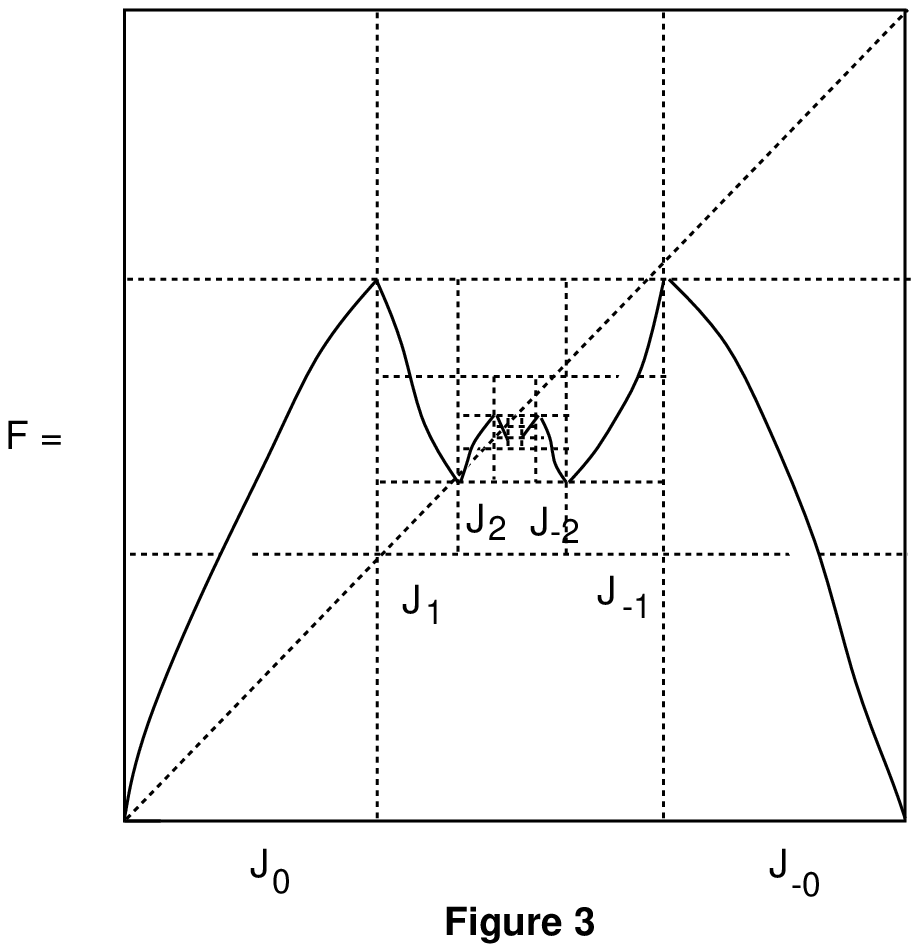}}

\vskip5pt
Let $g_{n}=(F|J_{n})^{-1}$ and $g_{-n}=(F|J_{-n})^{-1}$, $n=0$, $1$, $\cdots
$, be the inverse
branches of $F$. Suppose $w=i_{0}i_{1}\cdots i_{k-1} $ is a finite sequence of
${\bf Z}\cup \{ -0 \}$. We say it is admissible if the range $J_{i_{l}}$ of
$g_{i_{l}}$ is contained in the domain $F_{i_{l-1}}(J_{i_{l-1}})$ of
$g_{i_{l-1}}$ for $l=1$, $\cdots $, $k-1$.
For an admissible sequence $w = i_{0}i_{1}\cdots i_{k-1}$, we can define a
map 
$g_{w}=g_{i_{0}}\circ g_{i_{1}}\circ  \cdots \circ g_{i_{k-1}}$.  We use
$D(g_{w})$ to denote the domain of $g_{w}$ and use $|D(g_{w})|$ to denote
the length of the interval $D(g_{w})$. 

\vskip5pt
{\sc Definition 1.} {\em We say the induced Markov map $F$ has bounded distortion
property if there is a constant $C_{5}=C_{5}(f)>0$ such that 

\vskip5pt
(a) $|J_{n}|/ \cup_{l=n+1}^{\infty}|J_{l}|\geq C_{5}$ and 
$|J_{-n}|/ \cup_{l=n+1}^{\infty}|J_{-l}|\geq C_{5}$ for all the integers
$n\geq 0 $, and 

\vskip5pt
(b) $| \Big( N(g_{w})\Big) (x)|\leq C_{5}/|D(g_{w})|$ for
all $x$ in $D(g_{w})$ and all admissible $w$.} 

\vskip5pt
The reason we give this definition is the following lemma. First, we note that if
$g_{w}$ has the property (b) in Definition 1, then distortion of $g_{w}$ 
\[ \frac{|g_{w}(x)|}{|g_{w}(y)|} \leq \exp(C_{5})\]
for all $x$ and $y$ in $D(g_{w})$.

\vskip5pt
{\sc Lemma 5.} {\em Suppose $f$ and $g$ in ${\cal U}$ are two Feigenbaum
maps and $H$ is the conjugacy 
between $f$ and $g$. If both of the induced Markov maps $F$ and $G$ have 
bounded
distortion property, then $H$ is quasisymmetric.}

\vskip5pt
{\sl Proof.} It can be proved by almost the same arguments as that 
we used in the paper \cite{j2}. For
more details of the proof, the reader may refer to \cite{j2}.    

\vskip5pt
Suppose $g$ is a $C^{3}$-diffeomorphism 
with nonnegative Schwarzian derivative defined on $J$. 
We say $g$ has a {\em Koebe space} $C>0$ around $J$ if
$g$ can be extended to an interval $I\supset J$ as a $C^{3}$-diffeomorphism 
with nonnegative Schwarzian
derivative,  and moreover, 
the minimum of the ratios $|L|/|J|$ and $|R|/|J|$ is greater
than $C$ where $L$ and $R$ are the connected components of $I\setminus
J$.

\vskip5pt
{\sc Lemma 6.} {\em Suppose $g$ is a $C^{3}$-diffeomorphism with nonnegative Schwarzian
derivative defined on $J$. Moreover, suppose $g$ has a Koebe space $C>0$
around $J$. Then the nonlinearity $N(g)$ satisfies
\[ |\Big( N(g)\Big) (x)| \leq \frac{C'}{|J|}\]
where $C'=C'(C) >0$ is a constant.}

\vskip5pt
{\sl Proof.} It is a well-known lemma. We call it the $C^{3}$-Koebe
distortion lemma. See, for
example \cite[etc.]{bl,gj,ms,j1,su}. 

\vskip5pt
Now the proof of Theorem 1 concentrates on the next lemma.

\vskip5pt
{\sc Lemma 7.} {\em Suppose $f=h \circ Q_{t}$, for some $t>1$, in ${\cal U}$ is a Feigenbaum map  
and
$F$ is the Markov map induced by $f$.  Then $F$ has 
the bounded distortion property.}

\vskip5pt
{\sl Proof.} The condition (a) in Definition 1 is assured by Lemma 2. We only
need to check the condition (b) in Definition 1. 

\vskip5pt
Suppose $n$ is an integer $n\neq 0$, then $g_{n}$ can be extended to the interval $\Omega_{|n|} = L_{|n|}\cup
M_{|n|-1}$ as a $C^{3}$-diffeomorphism with nonnegative Schwarzian derivative 
(because $h$
is $C^{3}$-diffeomorphism with nonpositive Schwarzian derivative).  
For $g_{0}$ and $g_{-0}$,
without loss of generality, we may assume that they can be extended to the
interval $(-\infty , -1]\cup L_{0}$.  Lemma 1 and Lemma 2 now
assured that the map $g_{n}$ has a definite Koebe space $C_{7}>0$ around its domain
for all $n$ in ${\bf Z}\cup\{ -0\} $. We note that the intervals 
$\Omega_{|n|}$ are
nested for $|n|=0$, $1$, $\cdots$.

\vskip5pt
Suppose $w=i_{0}i_{1}\cdots
i_{k-1}$ is an admissible sequence of ${\bf Z}\cup \{ -0 \}$ and  
$g_{w} =g_{i_{0}}\circ g_{i_{1}}\circ \cdots \circ g_{i_{k-1}}$.
By the definition of an admissible sequence, one can check that (see Figure 1) 
\[ |i_{0}| \leq |i_{1}| \leq \cdots \leq |i_{k-1}|.\] 
Hence $g_{w}$ can be extended to the domain $\Omega_{|i_{k-1}|}$ as a
$C^{3}$-diffeomorphism with nonnegative Schwarzian derivative. Again the map
$g_{w}$ has the Koebe space $C_{7}$ around its domain. This implies that 
\[ |\Big( N(g_{w})\Big) (x)| \leq  \frac{C_{8}}{|D(g_{w})|}\]
for all $x$ in $D(g_{w})$ where $C_{8}=C_{8}(C_{7})> 0$ is a constant. 
This completes the proof of Lemma 7.

\vskip5pt
The arguments in Lemma 1 to Lemma 7 give the proof of Theorem 1.

\vskip20pt
Yunping Jiang

Institute for Mathematical Sciences 

SUNY at Stony Brook

Stony Brook, NY 11794

e-mail: jiang@math.sunysb.edu

\end{document}